\documentclass[a4paper,9pt] {article}
\usepackage[latin1]{inputenc}
\usepackage[T1]{fontenc}
\usepackage[english]{babel}
\usepackage{amsmath}
\usepackage{amsfonts}
\usepackage{amssymb}
\usepackage{graphicx}
\usepackage{color}
\definecolor{gris25}{gray}{0.55}
\usepackage{shadow}
\usepackage{makeidx}
\newcommand{\be}{\begin{equation}}
\newcommand{\ee}{\end{equation}}
\newcommand{\bd}{\begin{displaymath}}
\newcommand{\ed}{\end{displaymath}}

\newcommand{\ba}{\begin{eqnarray}}
\newcommand{\ea}{\end{eqnarray}}
\newcommand{\ban}{\begin{eqnarray*}}
\newcommand{\ean}{\end{eqnarray*}}

\newcommand{\R} {I\!\!R}
\newcommand{\E} {I\!\! E}
\newcommand{\N} {I\!\! N}

\newcommand{\ds}{\displaystyle}

\begin{document}
\title{Filtered derivative with p-value method for multiple change-points detection}

\date{ }
\maketitle \vspace{-1.5cm}
 \begin{center}
 Pierre, R.~BERTRAND${}^{1,2}$ and Mehdi FHIMA${}^{2}$\\
 ${}^{1}$  {\it INRIA Saclay, APIS Team} \\
${}^{2}$ {\it Laboratoire de Math\'ematiques, UMR CNRS 6620\\
\& University Clermont-Ferrand II, France}
\end{center}

\section*{Introduction}
In different applications (health, finance,...), abrupt changes on
the spectral density of long memory processes  provide relevant
information. In this work, we concern ourself with  off-line
detection. However, our method is close to the sliding window which
is typically a sequential analysis method.

We model data by  Gaussian processes with locally stationary
and long memory increments. By using a wavelet
analysis, 
one obtains a series with short
memory. We compare numerically the efficiency of
different methods for off-line detection of these changes, namely
penalized least square estimators introduced by Bai  and Perron
(1998) 
versus a
modification of the  filtered derivative introduced by Basseville
and Nikiforov (1993). 
The
enhancement consists in computing the p-value of every change point and
then apply an adaptive strategy.

Since estimation of abrupt changes on  spectral density  is a
specific problem, we first study a more standard model. In Section 1, we concern ourself to
off-line detection of abrupt  changes in the mean of independent
Gaussian variables with known variance and we numerically compare the efficiency
of the different estimators in this  case. In Section 2, we
recall the definition  of Gaussian  processes with locally stationary
increments and the properties of their wavelet coefficients. Then,
we compare the different off-line detection methods on simulated
locally fBm and present some results on real data.

\section{A toy model: off-line detection of abrupt change in the mean or independent Gaussian variables}
Let $\displaystyle (X_i)_{i=1,\dots,N}$ be a sequence of independent Gaussian
r.v. with mean $\mu_i$ and a known variance~$\sigma^2$. We assume that the
map $i\mapsto \mu_i$ is piecewise constant, {\em i.e.} there
exists a configuration $0=\tau_0<\tau_1<\dots<\tau_K<\tau_{K+1}=N$
such that $\mu_i= \mu_k$ for $\tau_k \le i < \tau_{k+1}$. The
integer $K$ corresponds to the number of changes. However, in any real life
situation, the number of abrupt changes $K$ is unknown, leading to
a problem of model selection.

There is a huge literature on this problem, see for instance the textbook of Basseville \& Nikiforov (1993).  Popular methods are those based on penalized least square criterion (PLSC). We refer to Birg\'e \& Massart (2006) for a good summary of the problem. Other classical references are  Lavielle \& Moulines (2000) or Lebarbier (2005) or Lavielle \& Teyssi\`ere (2006).

From a numerical point of view, the least square methods are based on  dynamic programming
algorithm. Thus we have to compute a matrix of size $N$. Therefore, the time and memory  complexity of these algorithms is in $O(N^2)$,
which
becomes an important limitation with the computer progress.
This has lead us to investigate the properties of a different algorithm.

\subsection*{Filtered derivative with p-value method (FDp-VM)}

Filtered derivative method is based on the difference between the
empirical mean computed on two sliding windows respectively at the right and
at the left of the index $k$, both of size $A$, see \cite{AH:94,
BN:93}. This difference corresponds to a sequence $(D(A,k))_{A
\leq k \leq N-A}$ defined by
$D(A,t)=\hat{\mu}(A,t)-\hat{\mu}(A,t-A)$ where
$\displaystyle \hat{\mu}(A,k)=\frac{1}{A} \sum_{j=k+1}^{k+A}X_j$ is the
empirical mean of $X$ on the (sliding) box $[k+1,\,k+A]$. These quantities can easily be
calculated by recurrence with complexity $O(N)$. It suffices to remark  that
$\displaystyle AD(A,k+1)=AD(A,k)+X_{k+A+1}-2X_k+X_{k-A+1}$.

From the other hand, note that
$(D(A,k))_{A\leq k \leq N-A}$
is  a  sequence of centered r.v., except in the vicinity of a change point $\tau_k$.
 In this case, there appears a "hat-function" of size $\delta_k:=\big(\mu_{k+1}-\mu_k\big)$
 approximatively located  between $\tau_k-A$ and $\tau_k+A$, see \cite{PRB:00}. 
Nevertheless, this sequence presents also  hats at
 points different from $(\tau_1,\ldots,\tau_K)$, see Figure
2 below. There are   false alarms.

In order to eliminate these false alarms, we propose to calculate
the  p-values  $\alpha_{k}$ associated to all detected change
points $\big(\tilde{\tau}_{k}\big)_{1\le k \le Kmax}$. Then, we
keep only the point corresponding to  a p-value lesser than a
critical level  $\alpha_{critic}$. For all $k\in [1, Kmax]$, the
p-values are calculated by using the formula $\alpha_k=\phi
\left(\frac{1}{\overline{\sigma}_{k}}\sqrt{\frac{A_k}{2}}|D(A_k,\tilde{\tau}_k)|\right)$
where $A_k$ is an adaptive window chosen as the minimum of the
distances  between $\tilde{\tau}_{k}$ and these two neighbors
$\tilde{\tau}_{k-1}$ and $\tilde{\tau}_{k+1}$, that is
$\displaystyle
A_k=|\tilde{\tau}_{k}-\tilde{\tau}_{k-1}|\wedge|\tilde{\tau}_{k+1}-\tilde{\tau}_{k}|$.
The function $\phi$ denotes the complementary cumulative
distribution function of the normal law and
$\overline{\sigma}_{k}$ denotes the empirical variance of $X$ on
the box $(\tilde{\tau}_{k-1}, \tilde{\tau}_{k+1})$.

Remark that $K_{max}$ is an integer fixed by the user. It represents the maximal number of
     change points. As soon as possible, $ K_{max}$ should be chosen  bigger than the true number of change points $K$.
     By convention, one sets $\tau_{0}=\tilde{\tau}_{0}=0$ and
     $\tau_{K+1}=\tilde{\tau}_{Kmax+1}=N$.

So, the novelty of this work consists in discriminating between
true and false alarms by attributing a p-value to each detected
change point. By keeping the change points with a p-value smaller
than $\alpha_{critic}$, one obtains the same precision as Lavielle \& Teyssi\`ere
(2006) or Lebarbier (2005).\\
Also, the main advantage of this method lies in its memory and
time complexity in $O(N)$.

\subsection*{A numerical simulation}
At first, we give an example on one sample. In the next
subsection, this example is plainly confirmed by Monte-Carlo
simulations. To begin with, for $N=5000$ we have simulated one
replication of a sequence of Gaussian random variable
$X_1,\ldots,X_N$ with variance $\sigma^2=1$ and mean
$\mu(i)=g(i/N)$ where $g$ is a piecewise-constant function with
five  change points such as $\delta_k \in [0.5,1.25]$, see
Figure~1 below.

\begin{figure}[htbp]
\begin{center}
\begin{tabular}{c}
\includegraphics[width=12cm,height=8cm]{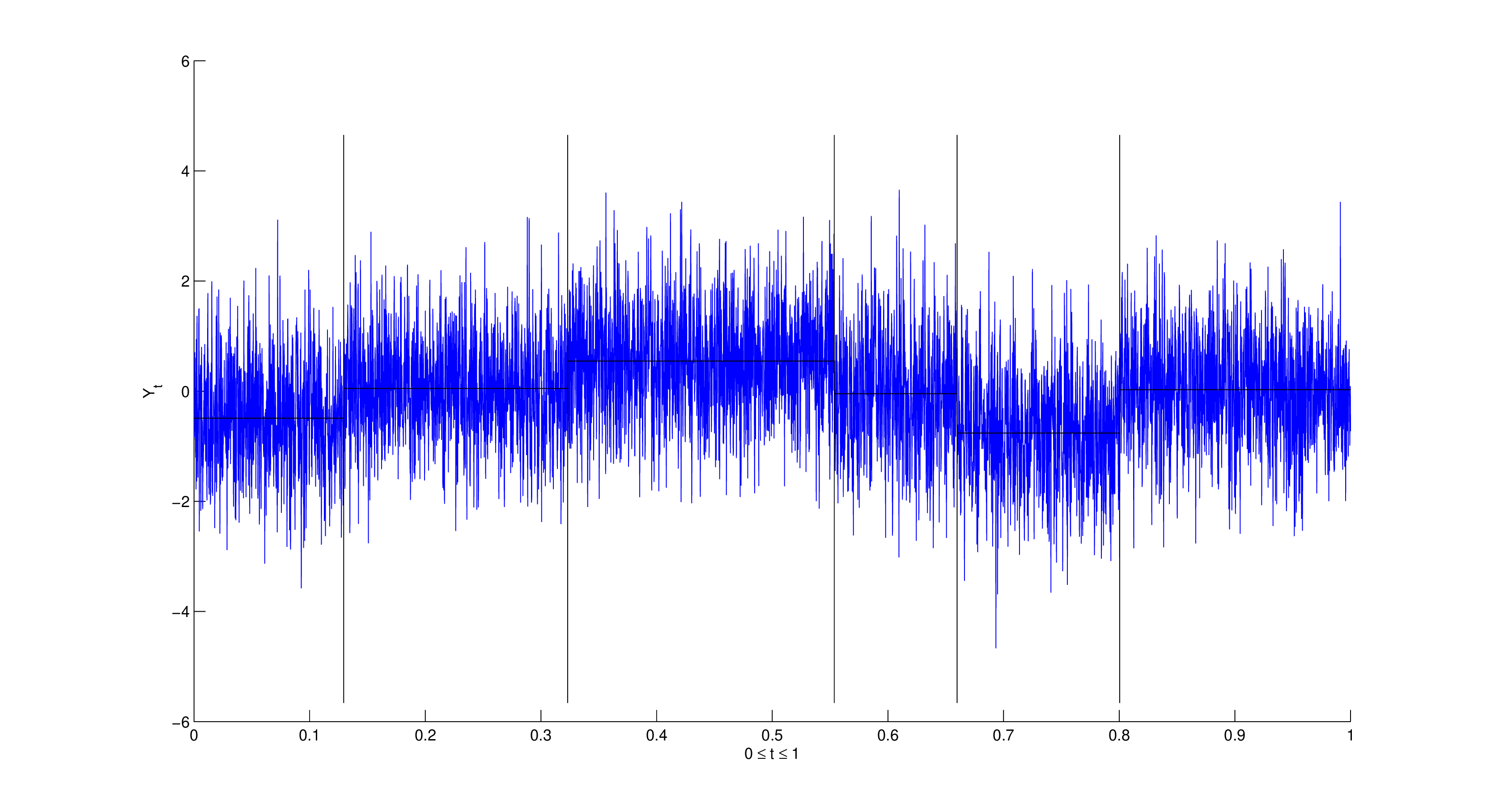}
\end{tabular}
\caption{\emph{The sequence $(X_i)_{0 \leq i \leq N}$ with change points represented by the vertical line and  means represented by the horizontal line
.}}
\end{center}
\end{figure}

\newpage

On this sample, we have computed the function $k\mapsto |D(A,k)|$
with $A=300$, see Figure~2.

\begin{figure}[htbp]
\begin{center}
\begin{tabular}{c}
\includegraphics[width=12cm,height=4cm]{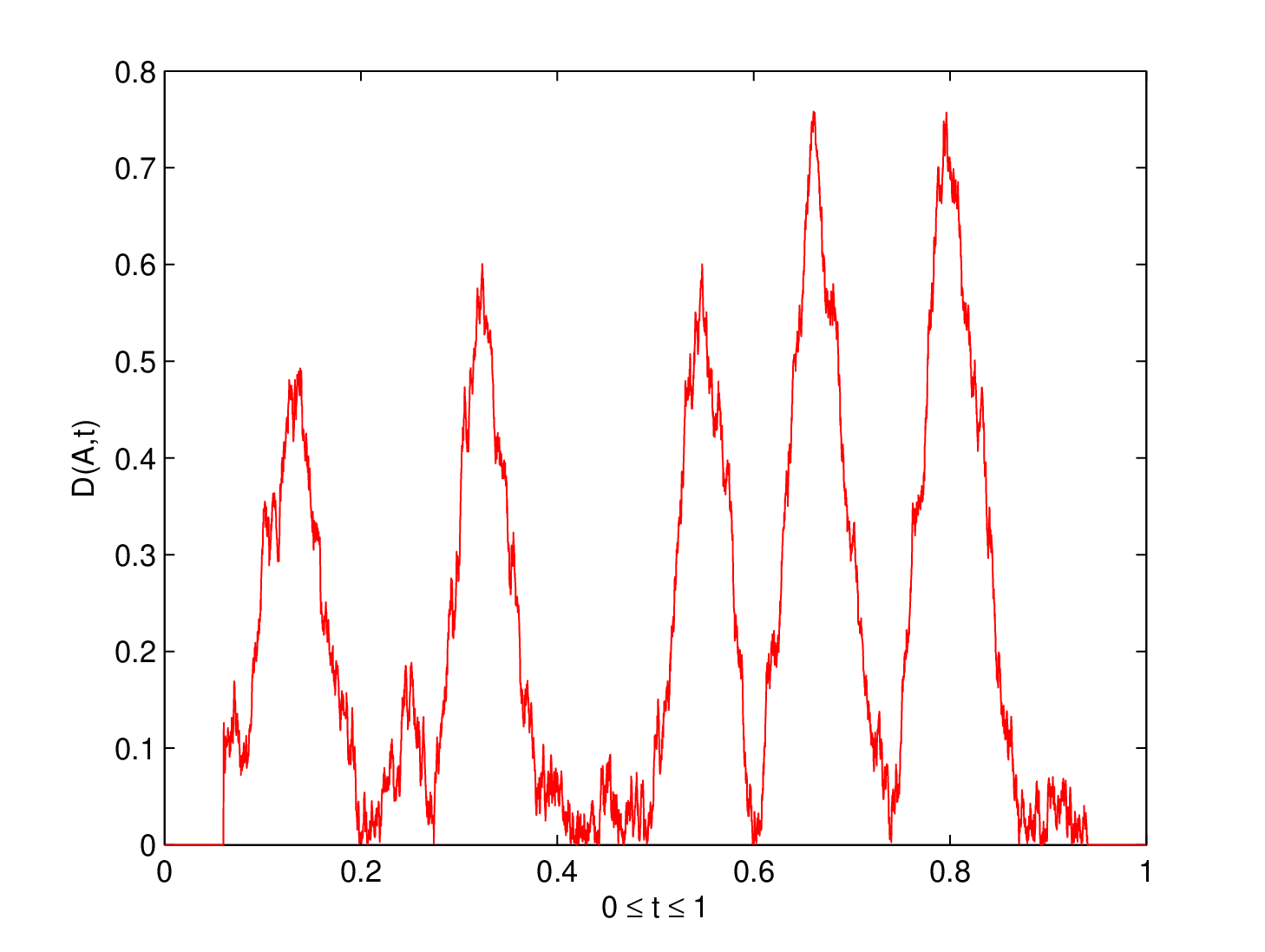}
\end{tabular}
\caption{\emph{The hat function}.}
\end{center}
\end{figure}

Both estimators penalized least square criterion (PLSC) 
and filtered derivative with p-value $\alpha_{critic}=10^{-4}$
provide good results, see Figure 3 and the Monte-Carlo simulation
below.\\

\begin{figure}[htbp]
\begin{center}
\begin{tabular}{c}
\includegraphics[width=12cm,height=4cm]{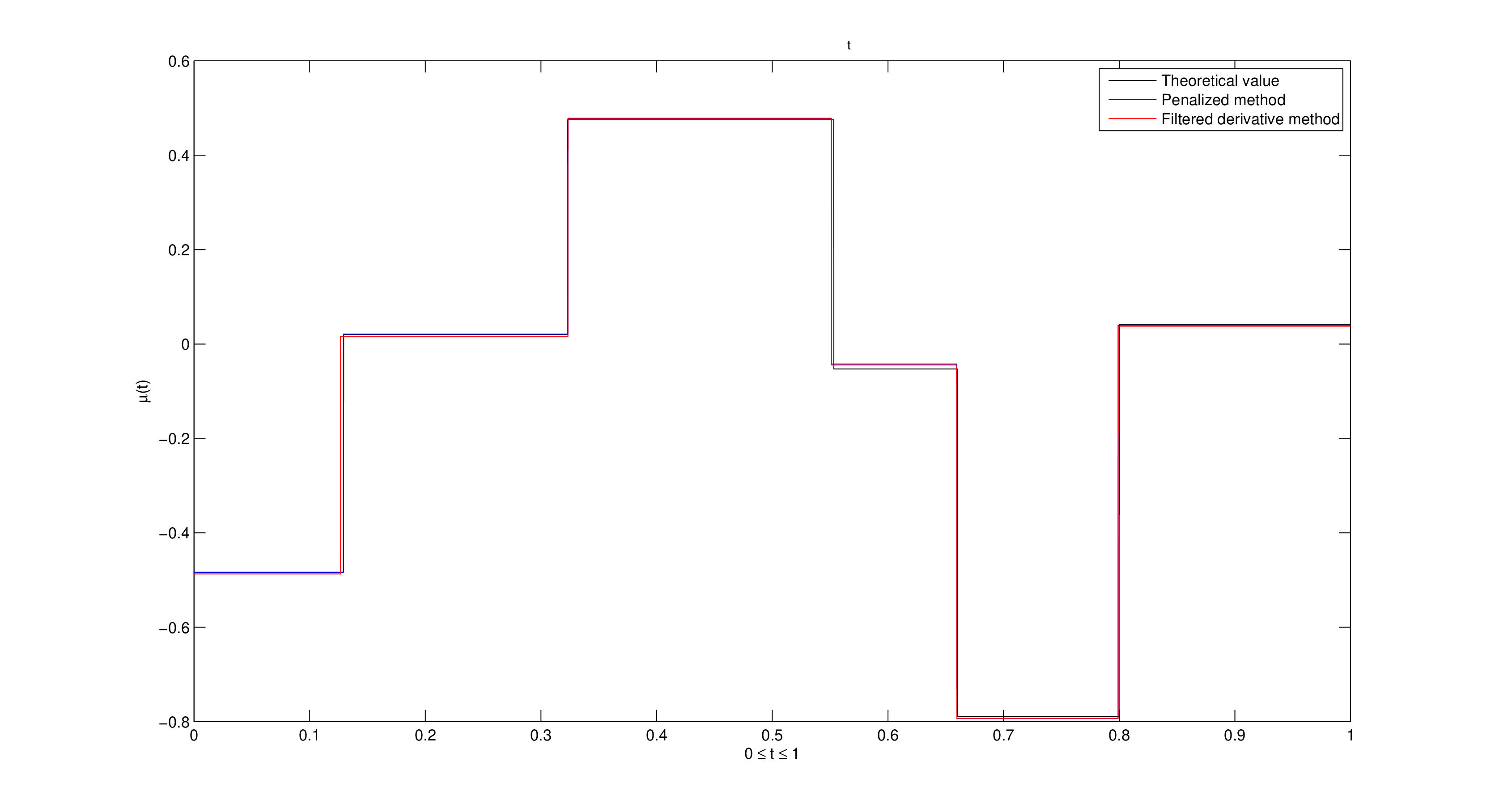}
\end{tabular}
\caption{\emph{Theoretical value of the piecewise-constant
function g (black), and its estimators given by PLSC method (blue)
and Filtered derivative method with p-value method (red)}.}
\end{center}
\end{figure}

\subsection*{Monte-Carlo simulation}
  In this subsection, we have made $M=1000$ simulations of independent copies of sequences
  of Gaussian r.v. $X_{0}^{(k)},\ldots,X_{N}^{(k)}$ with variance $\sigma^2=1$ and mean $\mu(i)=g(i/N)$, for $k=1,\dots,M$.
  On each sample, we apply the FDp-V method  and the PLSC method.
  We find the good number of changes in $98.1\%$ of all cases  for the first method and in $97.9\%$ for the second one.

\begin{figure}[htbp]
\begin{center}
\begin{tabular}{c}
\includegraphics[width=12cm,height=4cm]{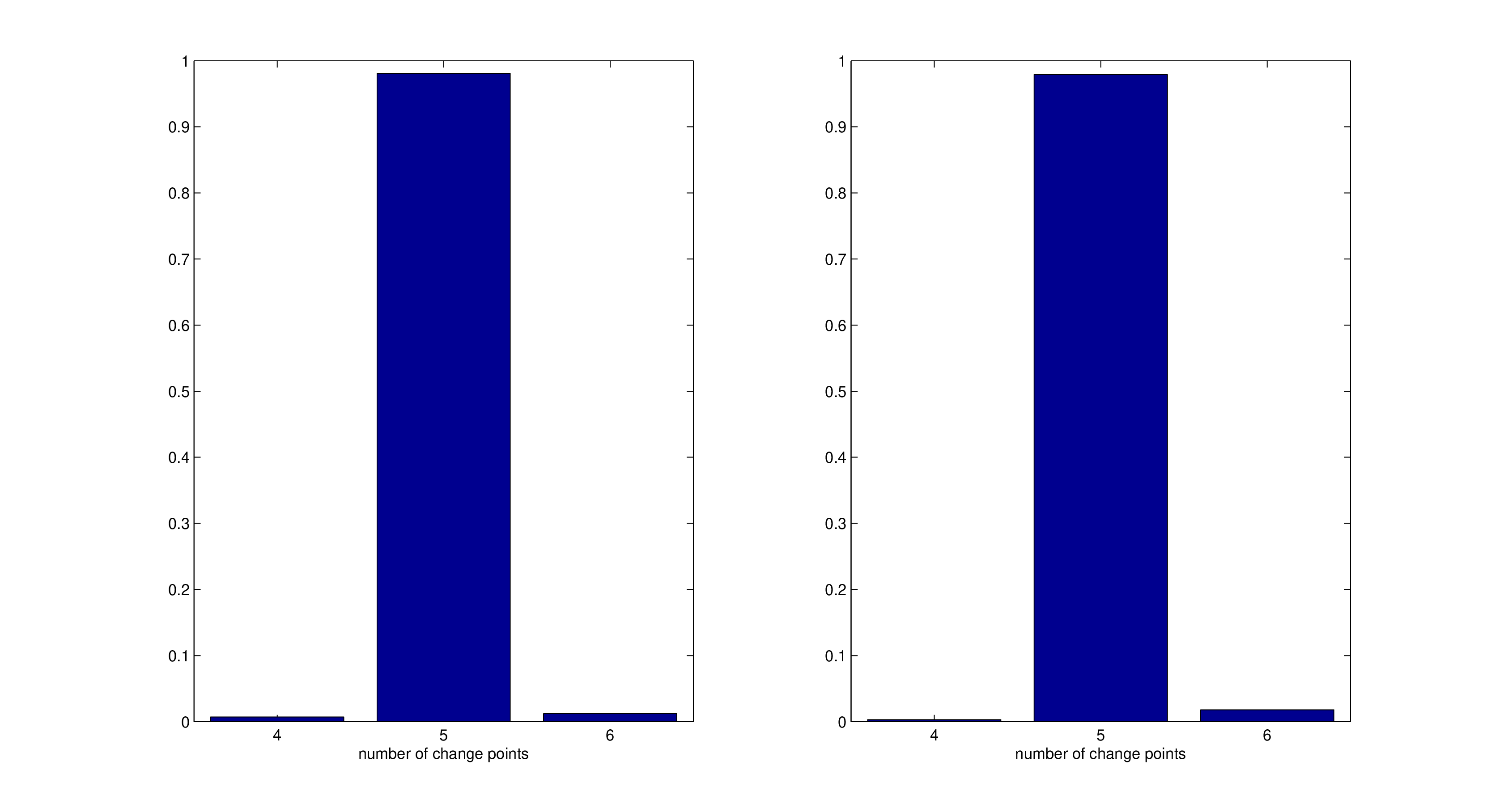}
\end{tabular}
\caption{\emph{Distribution of the estimated number of
change points $\widehat{K}$ for $M=1000$ realizations. Left:~Using
PLSC method. Right: Using Filtered derivative method}}
\end{center}
\end{figure}
\par\noindent
Then, we compute the mean errors. There are two kinds of mean error :
\begin{itemize}
  \item  Mean Integrate Square Error (MISE) defined as
$\displaystyle MISE=\E \left(\frac{1}{N+1} \sum_{i=0}^N
\left|\widehat{g}(i/N) - g(i/N)\right|^2 \right) = \E \|
\widehat{g} - g \|_{L^2([0,1])}^{2}$ which corresponds to the
$L^2([0,1])$ norm of the difference between the true function $g$
and the estimate function $\widehat{g}$. The estimate function is
obtained in two steps~: first we estimate the configuration of
change points $(\hat{\tau}_k)_{k=1,\dots, \hat{K}}$, then we
estimate the value of $\widehat{g}$ between two successive change
points as the empirical mean.
  \item Square Error on Change Points (SECP) defined as $\displaystyle SECP=\E  \left( \sum_{k=1}^K \left|\widehat{\tau}_k - \tau_k\right|^{2}\right)
$, in   the case where we have
found the good number of abrupt changes.
\end{itemize}
We have the following results by Monte Carlo simulation 
\begin{table}[htbp]
\begin{center}
\begin{tabular}{|c|c|c|}
  \hline
   \  & Square Error on Change Points   &  Mean Integrated Squared Error\\
  \hline
  FDp-V method & $1.1840 \times 10^{-4}$ & $0.0107$ \\
  \hline
  PLSC method & $1.2947 \times 10^{-4}$  & $0.0114$ \\
  \hline
\end{tabular}
\caption{Errors (SECP \& MISE) given by FDp-V method
and PLSC method}
\end{center}
\end{table}
\par
\noindent
Next, we compare the mean time complexity and the mean memory complexity.
We have written the two programs in Matlab and have runned it with computer
system which has the following characteristics: 1.8GHz processor
and $512$MB memory. The results concerning time and memory complexity
are given in Table 2.
\begin{table}[htbp]
\begin{center}
\begin{tabular}{|c|c|c|}
  \hline
   \  & Memory allocation (in Megabytes) &  CPU time (in second)\\
  \hline
  FDp-V method & $0.04$ MB & $0.005$ s \\
  \hline
  PLSC method & $200$ MB  & $240$ s \\
  \hline
\end{tabular}
\caption{Memory and time complexity of Filtered derivative method
and PLSC method}
\end{center}
\end{table}
\subsection*{A First conclusion}
On the one hand, both methods have the same accuracy in terms of
percentage of selection of the exact model, Square Error on the
configuration of change points or MISE. On the other hand, the
filtered derivative with p-value is less expensive in terms of time
complexity and memory complexity. Indeed, algorithm based on
Minimization of penalized least square criterion can  use $39\%$
of computer memory, while Filtered derivative method only  needs
$0.008\%$.
This plainly confirms the difference of time and
memory complexity, {\em i.e.} $O(N^2)$ {\em versus} $O(N)$.

Observe that algorithms based on penalized least square are considered by Davis et al. (2008) as maximizing {\em a posteriori} (MAP) criterion,
whereas filtered derivative is based on sliding window and could be adapted to sequential detection, see for instance Bertrand (2000) and Bertrand \& Fleury (2008).

\section{Segmentation on the spectral density estimation of some long memory processes}

\subsection*{Our model}
Let $X$ be a Gaussian centered process with stationary increments,
it is known, see Cram\'er \& Leadbetter  (1967), that this process
has the harmonizable representation $X(t)=\int_{\R}
\big(e^{it\xi}-1\big)f^{1/2}(\xi)dW(\xi)$ for all $t \in \R$ where $W(dx)$ is a
complex Brownian measure such that $X(t)$ is a real number for all $t\in \R$. This
process has long memory, but its wavelet coefficient $d_\psi(a,b)$
is a short memory Gaussian process where, for a scale and a shift
$(a,b)\in \R_+^*\times\R$ and a wavelet $\psi$ with a compact time
support $[L_1,L_2]$, one has defined \ba \label{def:coeff:wavelet}
d_\psi(a,b)&:=&  a^{-1/2} \int_{\R} \psi\left(
\frac{t-b}{a}\right)\,X(t)\,dt \ea Moreover, for any fixed scale
$a$,  $\ds b\mapsto d_\psi(a,b)$ is a stationary, centered,
Gaussian process  with  variance $\;\ds \mathcal{I}_{\psi}(a) :=
\int_{\R} |\widehat{\psi}(x)|^2\, f(x/a)\, dx$, see Bardet \&
Bertrand (2007). Next, we assume that the signal is a Gaussian
process, centered, with locally stationary increments given by the
representation formula \ba \label{repr:stat:harmonizable} X(t)=
\int_{\R} \big(e^{it\xi}-1\big)\, f^{1/2}(t,\xi) \,
dW(\xi),~~~~\mbox{for all}~~t \in \R, \ea where $\xi\mapsto
f(t,\xi)$ is an even and positive function, called spectral
    density piecewise constant, i.e., there exists a partition
$\tau_1 < \tau_2 < \dots < \tau_K$ and Hurst parameters
$H=(H_0,H_1,\ldots,H_K)$ such that $f(t,\xi) = f_k(\xi) =
C(H_k)|\xi|^{-2H_k-1}$ for $t\in [\tau_k, \tau_{k+1}[$ and
$C(H_k)=\pi^{-1}H_k\Gamma(2H_k)\sin(\pi H_k)$. Thus the series
$\log|d_\psi(a,b_i)|^2$ where $d_\psi(a,b)$ is defined by
(\ref{def:coeff:wavelet}) and $b_i=i\in \N$ has a piecewise
constant mean $\mu_i$  and  a finite known  variance, more
precisely, one has $\,\ds \log|d_\psi(a,b_i)|^2= \mu_i+ \zeta_i\,$
where $\zeta_i$ are weakly dependent r.v. of law $\ln |U|^2$ with
$U\sim\mathcal{N}(0,1)$ and $\,\ds \mu_i = \ln \int_{\R}
|\widehat{\psi}(x)|^2\, f_k(x/a)\,dx\;$ if $b_i\in
[\tau_k-aL_1,\tau_{k+1}-aL_2]$.

\subsection*{Numerical simulation}
First, for $T=10^5$, we have simulated one realization of process
$(X(t))_{t \in [0,T]}$ with five  change points
$\textbf{$\tau$}=\{12500,25496,43045,70083,82040\}$ and Hurst
parameters $H=(0.55,0.67,0.53,0.61,0.7,0.57)$. Let us stress that, after having changed the scale in order to obtain Hurst index belonging to $(0,1)$, the configuration of change points and means is the same as in Section 1.

\begin{figure}[htbp]
\begin{center}
\begin{tabular}{c}
\includegraphics[width=7cm,height=3cm]{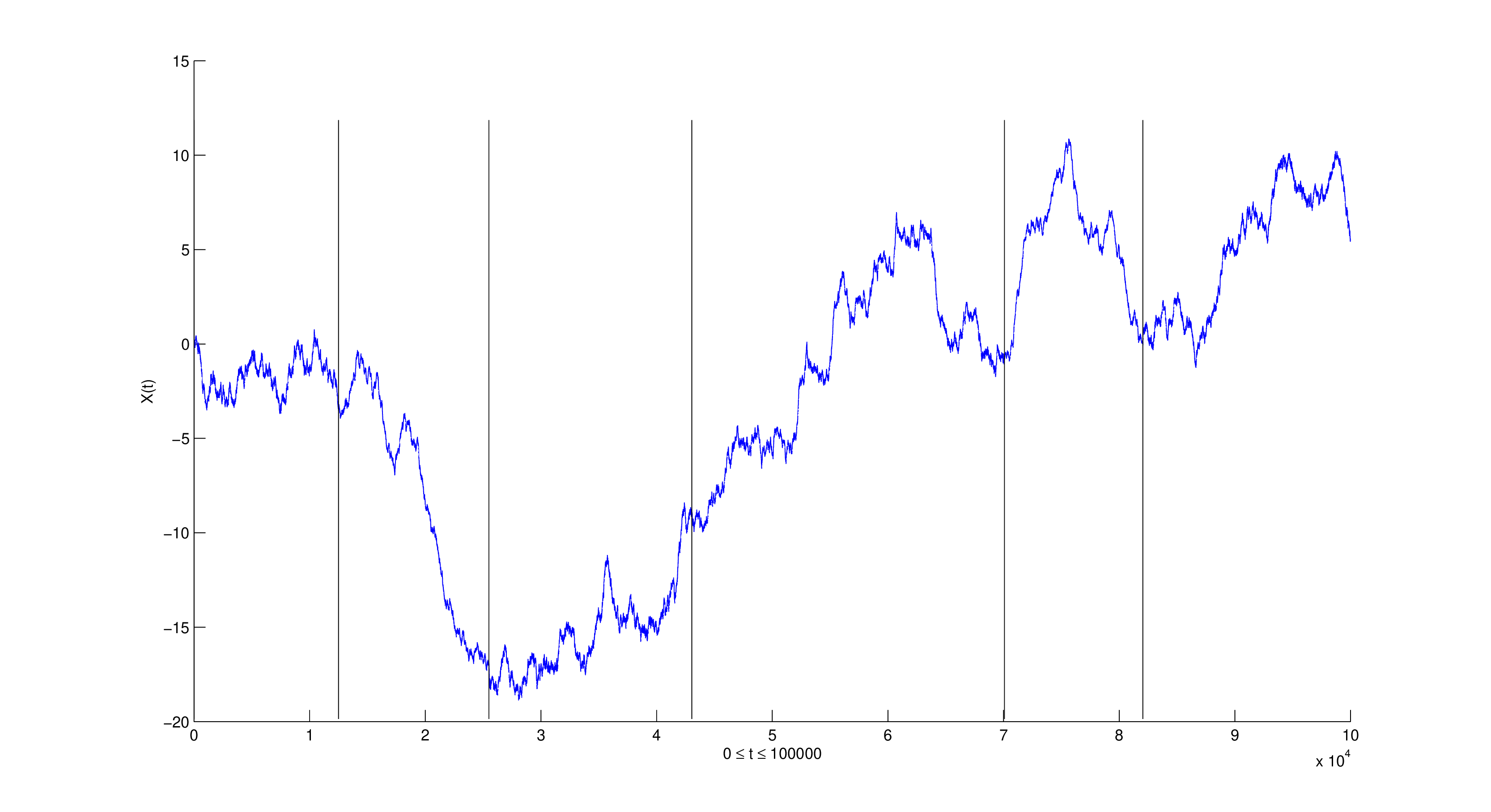}
\end{tabular}
\caption{\emph{One replication of the process $(X(t))_{t \in
[0,T]}$}}
\end{center}
\end{figure}

Next, for a frequency $1/a =0.2$Hz and by using the
Daubechies wavelet of order 6, we have calculated the wavelet
coefficients $(d_{\psi}(a,b_0),\ldots, d_{\psi}(a,b_N))$ with
$b_k=k$. Figure 5 below displays the sequence $(Y_0,\ldots,Y_N)$
where $Y_k=\log|d_{\psi}(a,b_k)|^2$. Then, we have calibrated the
Filtered derivative algorithm with $A=500$ and
$\alpha_{critic}=10^{-11}$. We observe that the detected change points
 perfectly fit the theoretical configuration of changes.
\begin{figure}[htbp]
\begin{center}
\begin{tabular}{c}
\includegraphics[width=7cm,height=3cm]{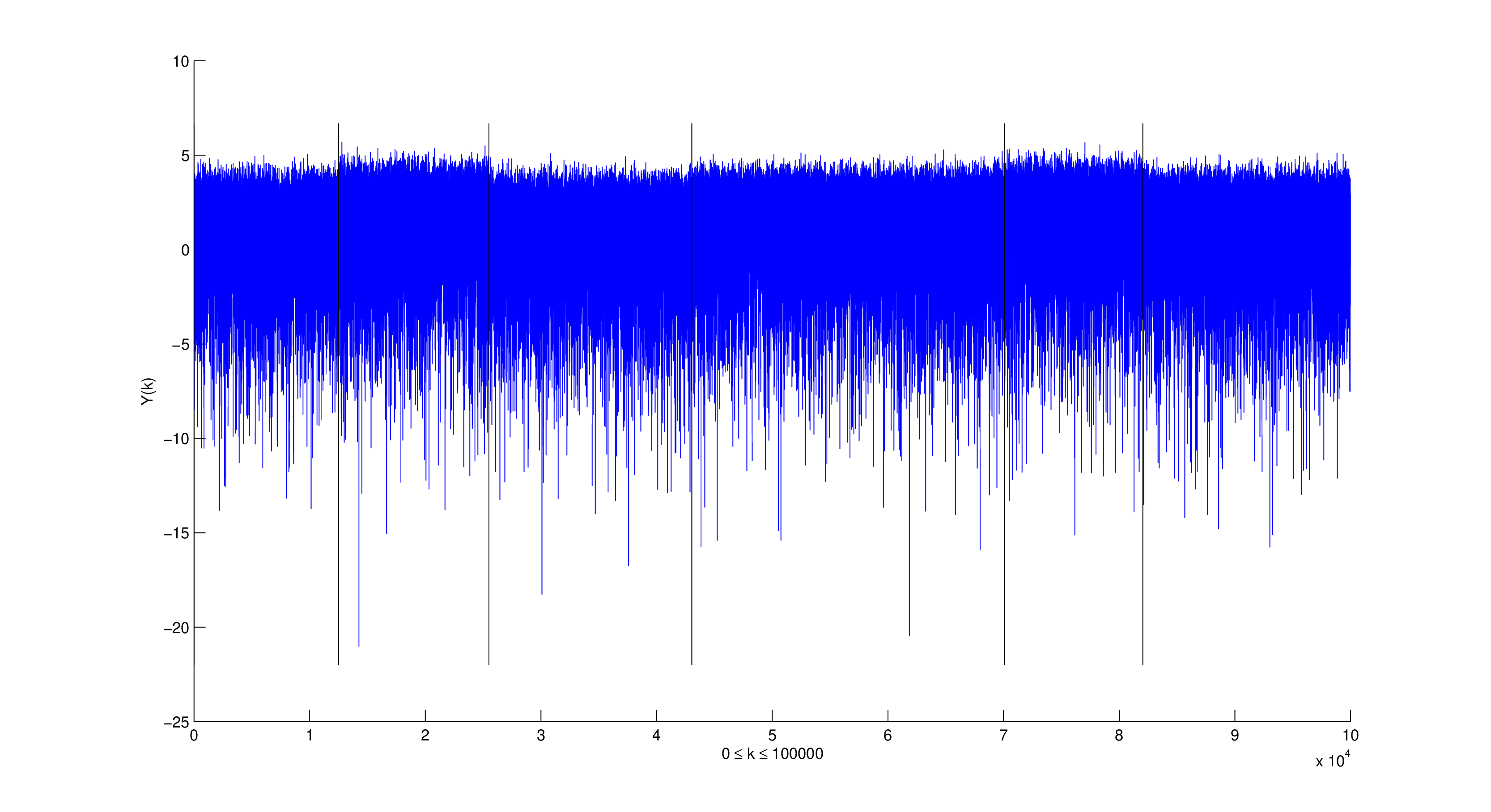}
\end{tabular}
\caption{\emph{Segmentation of the sequence $(Y_0,\ldots,Y_N)$}}
\end{center}
\end{figure}

Note that we can not  use penalized least square criterion due to the size of data, indeed  PLSC would have need $80$GB which is almost $150$ times our computer memory capacity.

\newpage

\section{Application to real data}
Recent measurement methods allow us to access to
electrocardiograms (ECG) for healthy people over a long period of
time: marathon runners, daily (24 hours) records, etc. These large
data sets allow us to characterize the variation of the heartbeat
rate in the parasympathetic  frequency band $(0.15\,Hz,0.5\,Hz)$.
According to the recommendations of the Task Force of
Cardiologists \cite{TaskForce:05}, this frequency band corresponds
to the parasympathetic system of control of the heartbeat.
Moreover, the spectral density of the heart beat time series
follows a power law, thus after having substract its mean, this
series  can be modelized by (2).
 Figure 7 provides an example of interbeat time series record on an healthy subject during 24 hours.   We have  calculated its
wavelet coefficients and used the Filtered derivative algorithm
with $A=500$ and $\alpha_{critic}=10^{-11}$. We obtain the
following
segmentation: $\tau=\{14435,21903,28003,31984,33377,37274,40470,42306,73153\}$

\begin{figure}[htbp]
\begin{center}
\begin{tabular}{c}
\includegraphics[width=7cm,height=3cm]{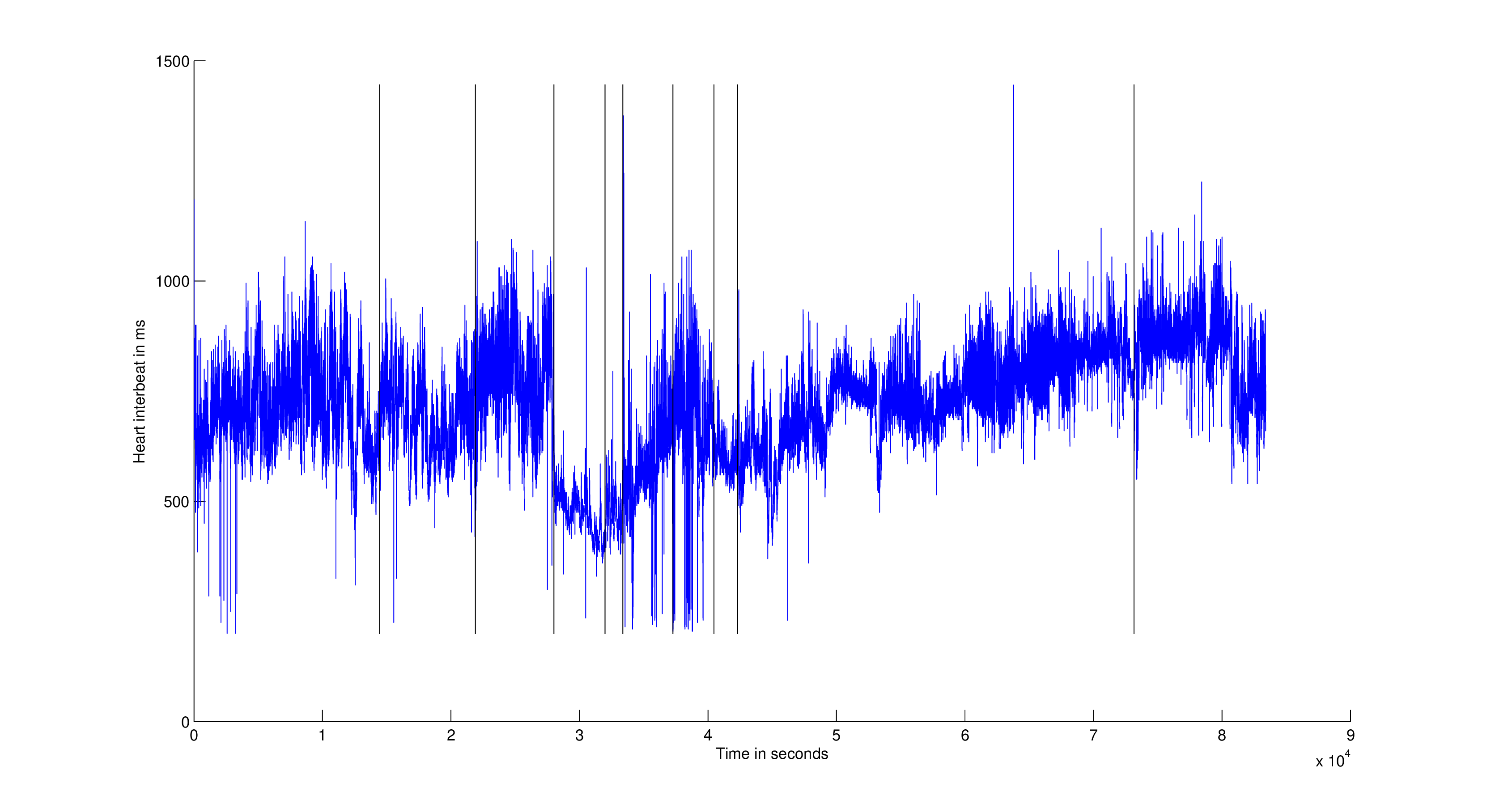}
\end{tabular}
\caption{\emph{Segmentation of the heart interbeat for healthy
subjects during a period of 24 hours}}
\end{center}
\end{figure}
\par\noindent
In future works, we will investigate sequential detection of change points of the Hurst index in connection with the
 cardiac behavior of sick subject.

\end{document}